\documentclass[11pt]{article}
\usepackage[table]{xcolor}
\usepackage{graphicx}
\usepackage{amsmath,amsfonts,amssymb}
\usepackage{graphicx}
\usepackage{multirow}
\usepackage{tikz,pgf}
\usepackage{url}
\usepackage{amsmath}
\usepackage{graphicx}
\usepackage{epsfig}
\usepackage{epstopdf}
\usepackage{color}
\usepackage{enumitem}
\usepackage{float}
\usepackage{algorithm}
\usepackage{chngcntr}
\counterwithin{figure}{section}
\def\disp{\displaystyle}

\def\tto{\;{\lower 1pt \hbox{$\rightarrow$}}\kern -10pt
\hbox{\raise 2pt \hbox{$\rightarrow$}}\;}
\def\h{\hfill\Box}

\def\Bar{\overline}
\def\ra{\rangle}
\def\la{\langle}

\def\epsilon{\varepsilon}
\def\B{\Bbb B}
\def\h{\hfill\Box}
\def\R{\Bbb R}

\def\N{\Bbb N}
\def\ox{\bar{x}}

\def\dom{\mbox{\rm dom}\,}

\def\h{\hfill\square}

\def\ph{\varphi}
\def\emp{\emptyset}

\def\oR{\Bar{\R}}

\def\ph{\varphi}
\def\emp{\emptyset}

\def\oR{\Bar{\R}}

\setlist[enumerate,1]{itemsep=0.0ex,parsep=0.5ex,label={\rm(\alph*)},leftmargin=*, align=left}
\newcounter{lk}

\usepackage[colorlinks=true]{hyperref}
\begin{document}
\begin{center}
{\sc\bf Solving  a Continuous Multifacility Location Problem by DC Algorithms}\\[2ex]
{\sc Anuj Bajaj}\footnote{ Department of Mathematics, Wayne State University, Detroit, Michigan 48202, USA (anuj.bajaj@wayne.edu ). Research of this author was partly supported by the USA National Science Foundation under grant DMS-1808978 and by the USA Air Force Office of Scientific Research grant \#15RT04.},
{\sc Boris S. Mordukhovich}\footnote{Department of Mathematics, Wayne State University, Detroit, Michigan 48202, USA (boris@math.wayne.edu). Research of this author was partly supported by the USA National Science Foundation under grants DMS-1512846 and DMS-1808978, by the USA Air Force Office of Scientific Research grant \#15RT04, and by Australian Research Council under grant DP-190100555.},
{\sc Nguyen Mau Nam}\footnote{Fariborz Maseeh Department of Mathematics and Statistics, Portland State University, Portland, OR 97207, USA (mnn3@pdx.edu). Research of this author was partly supported by the USA National Science Foundation under grant DMS-1716057.}, and {\sc Tuyen Tran}\footnote{Fariborz Maseeh Department of Mathematics and Statistics, Portland State University, Portland, OR 97207, USA (tuyen2@pdx.edu).}.\\[2ex]
\end{center}
\small{\bf Abstract.} The paper presents a new approach to solve multifacility location problems, which is based on mixed integer programming and algorithms for minimizing differences of convex (DC) functions. The main challenges for solving the multifacility location problems under consideration come from their intrinsic discrete, nonconvex, and nondifferentiable nature. We provide a reformulation of these problems as those of continuous optimization and then develop a new DC type algorithm for their solutions involving Nesterov's smoothing. The proposed algorithm is computationally implemented via MATLAB numerical tests on both artificial and real data sets. \\[1ex]

\newtheorem{Theorem}{Theorem}[section]
\newtheorem{Proposition}[Theorem]{Proposition}
\newtheorem{Remark}[Theorem]{Remark}
\newtheorem{Lemma}[Theorem]{Lemma}
\newtheorem{Corollary}[Theorem]{Corollary}
\newtheorem{Definition}[Theorem]{Definition}
\newtheorem{Example}[Theorem]{Example}
\renewcommand{\theequation}{\thesection.\arabic{equation}}
\normalsize

\section{Introduction and Problem Formulation}\label{intro}
\setcounter{equation}{0}
In the 17th century, Pierre de Fermat posed a problem of finding for a point that minimizes the sum of its Euclidean distances to three given points in the plane. The problem was soon solved by Evangelista Torricelli, and it is now known as the \emph{Fermat-Torricelli problem}. This problem and its extended version that involves a finite number of points in higher dimensions are examples of continuous {\em single facility} location problems. Over the years several generalized models of the Fermat-Torricelli type have been introduced and studied in the literature with practical applications to facility location decisions; see \cite{HM2015,Martini,n2,nars,nh,NamGiles} and the references therein. An important feature of single facility location problems and the problems studied in the aforementioned references is that only one \emph{center/server} has to be found to serve a finitely many \emph{demand points/customers}.

However, numerous practical applications lead to formulations of facility location problems in which more than one center must be found to serve a finite number of demand points. Such problems are referred to as \emph{multifacility location problems} (MFLPs). Given a finite number of demand points $a_1,\ldots,a_n$ in $\R^d$, we consider here the facility location in which $k$ centers $v_1,\ldots,v_k$ ($1\le k\le n$) in $\R^d$ need to be found to serve these demand points by assigning each of them to its nearest center and minimizing the total distances from the centers to the assigned demand points. In the case where $k=1$, this problem reduces to the generalized Fermat-Torricelli problem of finding a point that minimizes the sum of the distances to a finite number of given points in $\R^d$.

Let us formulate the problem under consideration in this paper as the following problem of \emph{mixed integer programming} with {\em nonsmooth} objective functions. It is convenient to use a variable $k\times d$-matrix $\mathbf V$ with $v_i$ as its $i$th row to store the centers to be found. We also use another variable $k\times n$-matrix $\mathbf U=[u_{i,j}]$ with $u_{i,j}\in\{0,1\}$ and $\sum_{i=1}^{k}u_{i,j}=1$ for $j=1,\ldots,n$ to assign demand points to the centers. The set of all such matrices is denoted by $\mathcal{U}$. Note that $u_{i,j}=1$ if the center $v_i$ is assigned to the demand point $a_j$ while $\sum_{i=1}^{k}u_{i,j}=1$ means that the demand point $a_j$ is assigned to only one center. Our goal is to solve the constrained optimization problem formulated as follows:
\begin{eqnarray}\label{main problem}
\begin{array}{ll}
&\mbox{\rm minimize }\mathcal{F}(\mathbf U,\mathbf V):=\sum_{i=1}^{k}\sum_{j=1}^{n}u_{i,j}^2\|a_{j}-v_{i}\|\\
&\mbox{\rm subject to }\mathbf U\in\mathcal{U}\;\;\mbox{\rm and }\;\mathbf V\in\R^{k\times d}.
\end{array}
\end{eqnarray}
Taking into account that $u_{i,j}\in\{0,1\}$, it is convenient to use $u^2_{i,j}$ instead of $u_{i,j}$ in the definition of the objective function $\mathcal{F}$; See Section 4.

Note that a similarly looking problem was considered by An, Minh and Tao \cite{ha} for different purposes. The main difference between our problem \eqref{main problem} and the one from \cite{ha} is that in \cite{ha} the squared Euclidean norm is used instead of the Euclidean norm in our formulation. From the point of applications this difference is significant; namely, using the Euclidean norm allows us to model the total distance in \emph{supply delivery}, while using the squared Euclidean norm is meaningful in {\em clustering}. Mathematically these two problems are essentially different as well. In addition to the challenging discrete nature and nonconvexity that both problems share, the objective function of our multifacility location problem \eqref{main problem} is {\em nondifferentiable} in contrast to \cite{ha}. This is yet another serious challenge from both theoretical and algorithmic viewpoints. Observe also that for $k=1$ our problem becomes the aforementioned generalized Fermat-Torricelli problem that does not have a closed-form solution, while the problem considered in \cite{ha} reduces to the standard problem of minimizing the sum of squares of the Euclidean distances to the demand points. The latter has a simple closed form solution given by the mean of the data points.

In this paper we develop the following algorithmic procedure to solve the formulated nonsmooth problem \eqref{main problem} of mixed integer programming:

{\bf(i)} Employ {\em Nesterov's smoothing} to approximate the nonsmooth objective function in \eqref{main problem} by a family of {\em smooth} functions, which are represented as {\em differences of convex} ones.

{\bf(ii)} Enclose the obtained smooth discrete problems into constrained problems of {\em continuous DC optimization} and then approximate them by unconstrained ones while using {\em penalties}.

{\bf(iii)} Solve the latter class of problems by developing an appropriate modification of the algorithm for minimizing differences of convex functions known as the {\em DCA}.

As a result of all the three steps above, we propose a {\em new algorithm} for solving the class of multifacility local problems of type \eqref{main problem}, verify its efficiency and implementation with MATLAB numerical tests on both artificial and real data sets.

Recall that the early developments on the DCA trace back to the work by Tao in 1986 with more recent results presented in \cite{AnNam}--\cite{Tao}, \cite{TA1}, \cite{TA2}, and the bibliographies therein. Nesterov's smoothing technique was introduced in his seminal paper \cite{n83} and was further developed and applied in many great publications; see, e.g., \cite{n,n18} for more details and references. The combination of these two important tools provides an effective way to deal with nonconvexity and nondifferentiability in many optimization problems encountered in facility location, machine learning, compressed sensing, and imaging. It is demonstrated in this paper in solving multifacility location problems of type \eqref{main problem}.

The rest of this paper is organized as follows. Section~\ref{Def} contains the basic definitions and some preliminaries, which are systematically employed in the text. In Section~\ref{Survey} we briefly overview two versions of the DCA, discuss their convergence, and present two examples that illustrate their performances.

Section~\ref{nest-smooth} is devoted to applying Nesterov's smoothing technique to the objective function of the multifacility location problem \eqref{main problem} and constructing in this way a smooth approximation of the original problem by a family of DC ones. Further, we reduce the latter smooth DC problems of discrete constrained optimization to unconstrained problems by using an appropriate penalty function method. Finally, the obtained discrete optimization problems are enclosed here into the DC framework of unconstrained continuous optimization.

In Section~\ref{TheAlgorithm} we proposed, based on the above developments, a new algorithm to solve the multifacility location problem \eqref{main problem} by applying the updated version of the DCA taken from Section~\ref{Survey} to the smooth DC problems of continuous optimization constructed in Section~\ref{nest-smooth}. The proposed algorithm is implemented in this section to solving several multifacility problems arising in practical modeling. Section~\ref{Conclusion} summarizes the obtained results and discusses some directions of future research.

\section{Basic Definitions and Preliminaries}\label{Def}
\setcounter{equation}{0}
For the reader's convenience, in this section we collect those basic definitions and preliminaries, which are largely used throughout the paper; see the books \cite{HUL,bmn,r} for more details and proofs of the presented results.

Consider the difference of two convex functions $g-h$ on a finite-dimensional space and assume that $g\colon\R^d\to\oR:=(-\infty,\infty]$ is extended-real-valued while $h\colon\R^d\to\R$ is real-valued on $\R^d$. Then a general problem of {\em DC optimization} is defined by:
\begin{equation}\label{DCP}
\mbox{\rm minimize }f(x):=g(x)-h(x),\quad x\in\R^d.
\end{equation}
Note that problem \eqref{DCP} is written in the unconstrained format, but---due to the allowed infinite value for $g$---it actually contains the domain constraint $x\in\dom(g):=\{u\in\R^d\;|\;g(u)<\infty\}$. Furthermore, the explicit constraints of the type $x\in\Omega$ given by a nonempty convex set $\Omega\subset\R^d$ can be incorporated into the format of \eqref{DCP} via the indicator function $\delta_\Omega(x)$ of $\Omega$, which equals $0$ for $x\in\Omega$ and $\infty$ otherwise. The representation $f=g-h$ is called a {\em DC decomposition} of $f$. Note that the class of DC functions is fairly large and include many nonconvex functions important in optimization. We refer the reader to the recent book \cite{m18} with the commentaries and bibliographies therein for various classes of nonconvex optimization problems that can be represented in the DC framework \eqref{DCP}.

Considering a nonempty (may not be convex) set $\Omega\subset\R^d$ and a point $x\in\R^d$, define the \emph{Euclidean projection} of $x$ to $\Omega$ by
\begin{equation}\label{proj}
P(x;\Omega):=\big\{w\in\Omega\;\big|\;\|x-w\|=d(x;\Omega)\big\},
\end{equation}
where $d(x;\Omega)$ stands for the \emph{Euclidean distance} from $x$ to $\Omega$, i.e.,
\begin{equation}\label{dist}
d(x;\Omega):=\inf\big\{\|x-w\|\;\big|\;w\in\Omega\big\}.
\end{equation}
Observe that $P(x;\Omega)\ne\emp$ for closed sets $\Omega$ while being always a singleton if the set $\Omega$ is convex.

Given further an extended-real-valued and generally nonconvex function $\ph\colon\R^d\to\oR$, the {\em Fenchel conjugate} of $\ph$ is defined by
\begin{equation*}
\ph^*(v):=\sup\big\{\la v,x\ra-\ph(x)\;\big|\;x\in\R^d\big\},\quad v\in\R^d.
\end{equation*}
If $\ph$ is {\em proper}, i.e., $\dom(\ph)\ne\emp$, its Fenchel conjugate $\ph^*\colon\R^d\to\oR$ is automatically convex.

The {\em subdifferential} of $\ph\colon\R^d\to\oR$ at $\ox\in\dom(\ph)$ is the set of subgradients given by
\begin{equation}\label{sub}
\partial\ph(\ox):=\big\{v\in\R^d\;\big|\;\la v,x-\ox\ra\le\ph(x)-\ph(\ox)\;\;\mbox{\rm for all }\;x\in\R^d\big\}.
\end{equation}
If $\ox\notin\dom(\ph)$, we let $\partial\ph(\ox):=\emp$. Recall that for functions $\ph$ differentiable at $\ox$ with the gradient $\nabla\ph(\ox)$ we have $\partial\ph(\ox)=\{\nabla\ph(\ox)\}$. \vspace*{0.05in}

The following proposition gives us a two-sided relationship between the Fenchel conjugates and subgradients of convex functions.

\begin{Proposition}\label{characterization} Let $\ph\colon\R^d\to\oR$ be a proper, lower semicontinuous, and convex function. Then $v\in\partial\ph^*(y)$ if and only if
\begin{equation}\label{e1}
v\in\mbox{\rm argmin}\,\big\{\ph(x)-\la y,x\ra\;\big|\;x\in\R^d\big\}.
\end{equation}
We have furthermore that $w\in\partial\ph(x)$ if and only if
\begin{equation}\label{e2}
w\in\mbox{\rm argmin}\,\big\{\ph^*(y)-\la x,y\ra\;\big|\;y\in\R^d\big\}.
\end{equation}
\end{Proposition}
{\bf Proof.} To verify the first assertion, suppose that \eqref{e1} is satisfied and then get $0\in \partial\psi(v)$, where $\psi(x):=\ph(x)-\la y,x\ra$ as $x\in\R^d$. It tells us that
\begin{equation*}
0\in\partial\ph(v)-y,
\end{equation*}
and hence $y\in\partial\ph(v)$, which is equivalent to $v\in\partial\ph^*(y)$ due to the biconjugate relationship $\ph^{**}=\ph$ valued under the assumptions made.

In the opposite way, assuming $v\in\partial\ph^*(y)$ gives us by the proof above that $0\in\partial\psi(v)$, which clearly yields \eqref{e1} and thus justifies the first assertion.

To verify the second assertion, suppose that \eqref{e2} holds and then get $0\in\partial\psi(w)$, where $\psi(y):=\ph^*(y)-\la x,y\ra$ as $y\in\R^d$. This clearly implies that
\begin{equation*}
0\in\partial\ph^*(w)-x,
\end{equation*}
and hence $x\in\partial\ph^*(w)$, which is equivalent to $w\in\partial\ph(x)$ due to the biconjugate relationship. The proof of the opposite implication in \eqref{e2} is similar to the one given above. $\h$

Finally in this section, recall that for a given variable matrix $\mathbf U\in\R^{k\times n}$ as in the optimization problem \eqref{main problem}, the \emph{Frobenius norm} on $\R^{k\times n}$ is defined by
\begin{equation}\label{frob}
\|\mathbf U\|_{F}:=\sqrt{\sum_{i=1}^{k}\sum_{j=1}^{n}|u_{i,j}|^2}.
\end{equation}

\section{Overview of the DCA and Some Examples}\label{Survey}
\setcounter{equation}{0}
In this section we first briefly overview two algorithms of the DCA type to solve DC problems \eqref{DCP} while referring the reader to \cite{TA1,TA2} for more details and further developments. Then we present numerical examples illustrating both algorithms.\\[1ex]
{\bf Algorithm~1: DCA-1}.
\begin{center}
\begin{tabular}{| l |}
\hline
{\small INPUT}: $x_0\in\R^d$, $N\in\N$.\\
{\bf for} $l=1,\ldots,N$ {\bf do}\\
\qquad\small{Find} $y_{l-1}\in \partial h(x_{l-1})$.\\\textcolor[rgb]{1.00,0.50,0.50}{}
\qquad\small{Find} $x_{l}\in\partial g^*(y_{l-1})$.\\
{\bf end for}\\
{\small OUTPUT}: $x_{N}$.\\
\hline
\end{tabular}
\end{center}

Since the convex function $h\colon\R^d\to\R$ in \eqref{DCP} is real-valued on the whole space $\R^d$, we always have $\partial h(x)\ne\emp$ for all $x\in\R^d$. At the same time, the other convex function $g\colon\R^d\to\oR$ in \eqref{DCP} is generally extended-real-valued, and so the subdifferential of its conjugate $g^*$ may be empty. Let us present an efficient condition that excludes this possibility.
Recall that a function $g\colon\R^d\to\oR$ is \emph{coercive} if
\begin{equation*}
\lim_{\|x\|\to\infty}\frac{g(x)}{\|x\|}=\infty.
\end{equation*}

\begin{Proposition}\label{nonemp-sub} Let $g\colon\R^d\to\oR$ be a proper, lower semicontinuous, and convex function. If in addition $g$ is coercive, then $\partial g^*(v)\ne\emp$ for all $v\in\R^d$.
\end{Proposition}
{\bf Proof.} Since $g$ is proper, the conjugate function $g^*$ takes values in $(-\infty,\infty]$ being convex on $\R^d$. Taking into account that $g$ is also and lower semicontinuous and invoking the aforementioned biconjugate relationship, we find $w\in\R^d$ and $c\in\R$ such that
\begin{equation}\label{sub-est}
c+\la w,x\ra\le g(x)\;\;\mbox{\rm for all }\;x\in\R^d.
\end{equation}
The coercivity property of $g$ ensures the existence of $\eta>0$ for which
\begin{equation*}
\|x\|\big(\|w\|+1\big)\le g(x)\;\;\mbox{\rm whenever }\;\|x\|\ge\eta.
\end{equation*}
It follows furthermore that
\begin{equation*}
\sup\big\{\la v,x\ra-g(x)\;\big|\;\|x\|\ge\eta\big\}\le-\|x\|\;\;\mbox{\rm for any }\;v\in\R^d.
\end{equation*}
By using \eqref{sub-est}, we arrive at the estimates
\begin{equation*}
\sup\big\{\la v,x\ra-g(x)\;\big|\;\|x\|\le\eta\big\}\le\sup\big\{\la v,x\ra-\la w,x\ra-c\;\big|\;\|x\|\le\eta\big\}<\infty,\quad v\in\R^d.
\end{equation*}
This tells us that $g^*(v)<\infty$, and therefore $\dom(g^*)=\R^d$. Since $g^*$ is a convex function with finite values, it is continuous on $\R^d$ and hence $\partial g^*(v)\ne\emp$ for all $v\in\R^d$. $\h$

To proceed further, recall that a function $\ph\colon\R^d\to\oR$ is $\gamma$-\emph{convex} with a given modulus $\gamma\ge 0$ if the function $\psi(x):=\ph(x)-\frac{\gamma}{2}\|x\|^2$ as $x\in\R^d$ is convex on $\R^d$. If there exists $\gamma>0$ such that $\ph$ is $\gamma-$convex, then $\ph$ is called \emph{strongly convex} on $\R^d$.

We also recall that a vector $\ox\in\R^d$ is a \emph{stationary point} of the DC function $f$ from \eqref{DCP} if
\begin{equation*}
\partial g(\ox)\cap\partial h(\ox)\ne\emp.
\end{equation*}
The next result, which can be derived from \cite{TA1,TA2}, summarizes some convergence results of the DCA. Deeper studies of the convergence of this algorithm and its generalizations involving the  Kurdyka-Lojasiewicz (KL) inequality are given in \cite{AnNam,TAN}.

\begin{Theorem}\label{dsa-conv} Let $f$ be a DC function taken from \eqref{DCP}, and let $\{x_l\}$ be an iterative sequence generated by Algorithm~{\rm 1}. The following assertions hold:
\begin{enumerate}
\item The sequence $\{f(x_l)\}$ is always monotone decreasing.
\item Suppose that $f$ is bounded from below, that $g$ is lower semicontinuous and $\gamma_1$-convex, and that $h$ is $\gamma_2$-convex with $\gamma_1+\gamma_2>0$. If $\{x_l\}$ is bounded, then the
limit of any convergent subsequence of $\{x_l\}$ is a stationary point of $f$.
\end{enumerate}
\end{Theorem}

In many practical applications of Algorithm~1, for a given DC decomposition of $f$ it is possible to find subgradient vectors from $\partial h(x_l)$ based on available formulas and calculus rules of convex analysis. However, it may not be possible to explicitly calculate an element of $\partial g^*(y_l)$. Such a situation requires either constructing a more suitable DC decomposition of $f$, or finding $x_{l+1}\in\partial g^*(y_l)$  by using the description of Proposition~\ref{characterization}. This leads us to the following modified version of the DCA.

{\bf Algorithm~2: DCA-2}.
\begin{center}
\begin{tabular}{| l |}
\hline
{\small INPUT}: $x_0\in\R^d$,\;$N\in\N$\\
{\bf for}\;$l=1,\ldots,N$\;{\bf do}\\[1ex]
\qquad Find $y_{l-1}\in \partial h(x_{l-1})$\\[1ex]
\qquad Find $x_{l}$ by solving the problem:\\[1ex]
\qquad\qquad\qquad\qquad$\mbox{\rm minimize}\;\ph_l(x):=g(x)-\la y_{l-1},x\ra,\;x\in\R^d.$\\[1ex]
{\bf end for}\\
{\small OUTPUT}: $x_{N}$\\
\hline
\end{tabular}
\end{center}\vspace*{0.05in}

Let us now present two examples illustrating the performances of Algorithms~1 and 2. The first example concerns a polynomial function of one variable.

\begin{Example}\label{1D} {\rm Consider the function $f\colon\R\to\R$ given
\begin{equation*}
f(x):=x^4-2x^2+2x-3\;\;\mbox{\rm for }\;x\in\R.
\end{equation*}
This function admits the DC representation $f=g-h$ with $g(x):=x^4$ and $h(x):=2x^2-2x+3$. To minimize $f$, apply first the {\em gradient method} with constant stepsize. Calculating the derivative of $f$ is $f^\prime(x)=4x^3-4x+2$ and picking any starting point $x_0\in\R$, we get the sequence of iterates
\begin{equation*}
x_{l+1}=x_l-t(4x_l^3-4x_l+2)\;\;\mbox{\rm for }\;l=0,1,\ldots.
\end{equation*}
constructed by the gradient method with stepsize $t>0$. The usage of the DC {\em Algorithm}~1 (DCA-1) gives us $y_l=\nabla h(x_l)=4x_l-2$ and then
$g^*(x)=3(x/4)^{4/3}$ with $\nabla g^*(x)=(x/4)^{1/3}$. Thus the iterates of DCA-1 are as follows:
\begin{equation*}
x_{l+1}=\nabla g^*(y_l)=\Big(\frac{y_l}{4}\Big)^{1/3}=\Big(\frac{4x_l-2}{4}\Big)^{1/3}=\Big(\frac{2x_l-1}{2}\Big)^{1/3},\quad l=0,1,\ldots.
\end{equation*}
\begin{figure}[H]
\centering
\includegraphics[scale=0.5]{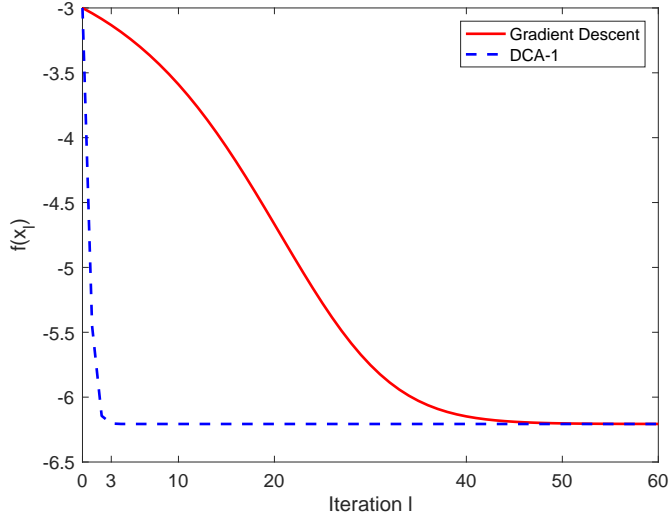}
\caption{Convergence of the gradient method and DCA-1.}
\label{convergencegraph}
\end{figure}
Figure~\ref{convergencegraph} provides the visualization and comparison between the DCA-1 and the gradient method. It shows that for $x_0=0$ and $t=0.01$ the DCA-1 exhibits much faster convergence.}
\end{Example}

The next two-dimensional example illustrates the performance of the DCA-2.

\begin{figure}[H]
\centering
\begin{minipage}{.55\textwidth}
\centering
\includegraphics[scale=0.53]{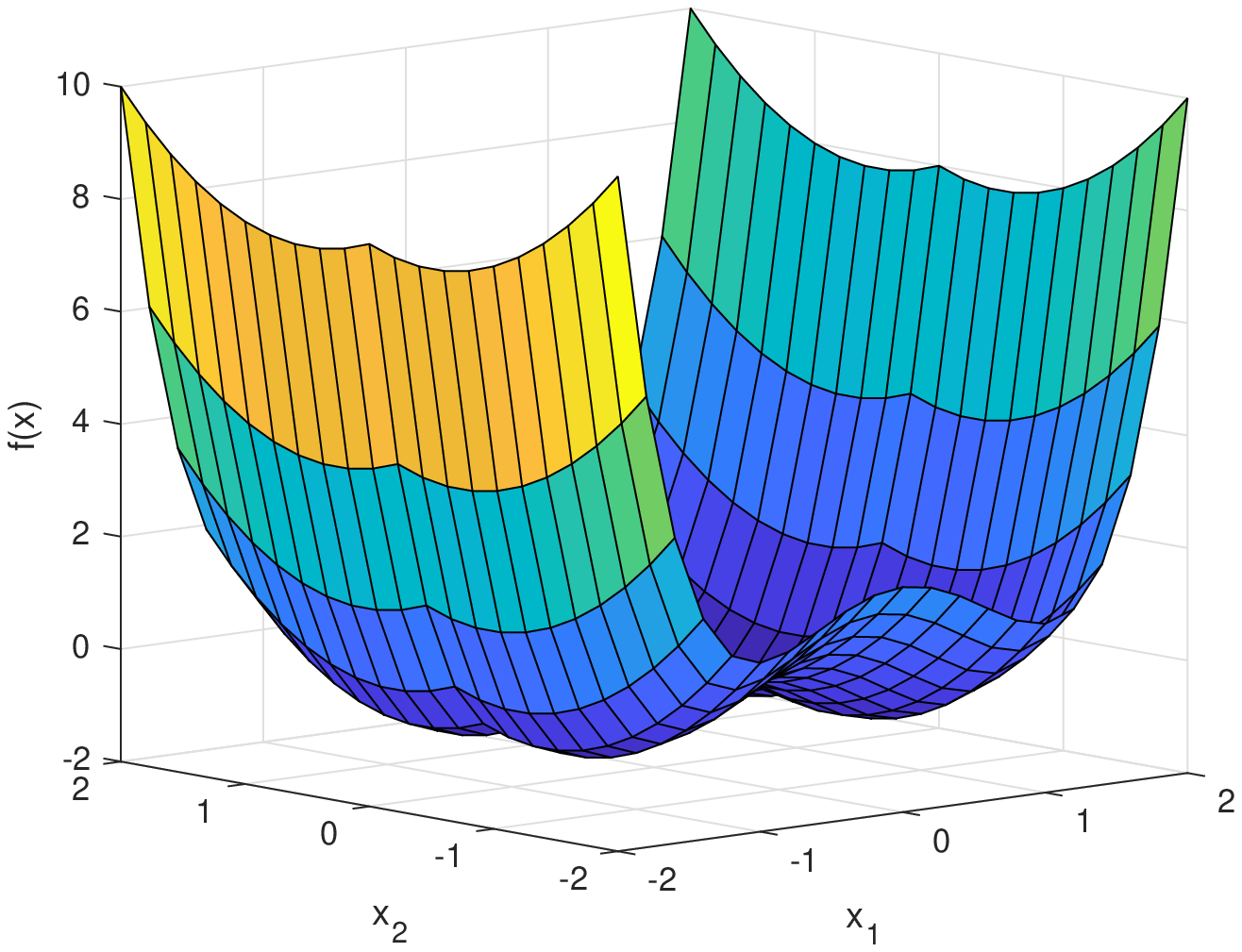}
\caption{Function graph for Example~\ref{2D}.}
\label{graph3D}
\end{minipage}
\begin{minipage}{.5\textwidth}
\centering
\includegraphics[scale=0.4]{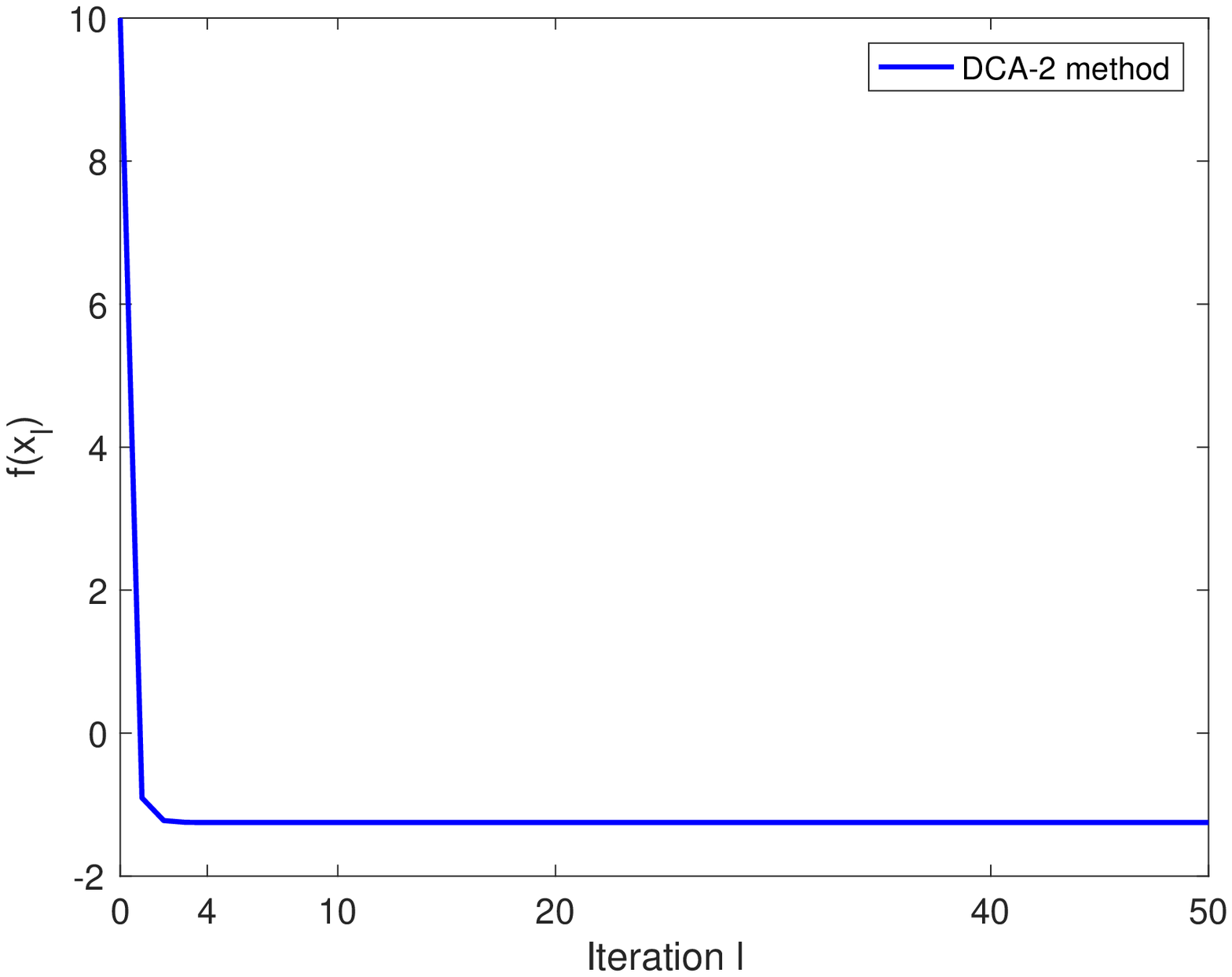}
\caption{Convergence of the DCA-2.}
\label{convergencegraph2D}
\end{minipage}
\end{figure}

\begin{Example}\label{2D} {\rm Consider the nonsmooth optimization problem defined by
\begin{equation*}
\mbox{\rm minimize }\;f(x_1,x_2):=x_1^4+x_2^2-2x_1^2-|x_2|\;\;\mbox{\rm over }\;x=(x_1,x_2)\in\R^2.
\end{equation*}
The graph of the function $f$ is depicted in Figure~\ref{graph3D}. Observe that this function has four global minimizers, which are $(1,0.5),(1,-0.5),(-1,0.5)$, and $(-1,-0.5)$. It is easy to see that $f$ admits a DC representation $f=g-h$ with $g(x_1,x_2):=x_1^4+x_2^2$ and $h(x_1,x_2):=2x_1^2+|x_2|$. We get the gradient
$\nabla g(x)=[4x_1^3,2x_2]^T$ and the Hessian
\begin{center}
$\nabla^2g(x)=
\begin{bmatrix}
12{x_1}^2&0\\
0&2
\end{bmatrix}$,
\end{center}
while an explicit formula to calculate $\partial g^*(y_l)$ is not available. Let us apply the {\rm DCA}-2 to solve this problem. The subdifferential of $h$ is calculated by
\begin{equation*}
\partial h(x)=\big[4x_1,\mbox{\rm sign}(x_2)\big]^T\;\;\mbox{\rm for any }\;x=(x_1,x_2)\in\R^2.
\end{equation*}
Having $y_l$, we proceed with solving the subproblem
\begin{equation}\label{subpr}
\mbox{\rm minimize }\;\ph_l(x):=g(x)-\la y_l,x\ra\;\;\mbox{\rm over }\;x\in\R^2
\end{equation}
by the classical Newton method with $\nabla^2\ph_l(x)=\nabla^2g(x)$ and observe that the DCA-2 shows its superiority in convergence with different choices of initial points. Figure~\ref{convergencegraph2D} presents the results of computation by using the DCA-2 with the starting point $x_0=(-2,2)$ and employing the Newton method with $\epsilon=10^{-8}$ to solve subproblem \eqref{subpr}.}
\end{Example}

\section{Smooth Approximation by Continuous DC Problems}\label{nest-smooth}
\setcounter{equation}{0}
In this section we first employ and further develop Nesterov's smoothing technique for the case of multifacility location problem \eqref{main problem}. Then we enclose the family of DC mixed integer programs obtained in this way into a class of smooth DC problems of continuous optimization. The suggested procedures are efficiently justified by deriving numerical estimates expressed entirely via the given data of the original problem \eqref{main problem}.

We begin with recalling the following useful result on Nesterov's smoothing related to the problem under consideration, which is taken from \cite[Proposition~2.2]{Nam2018}.

\begin{Proposition}\label{Prop3.1} Given any $a\in\R^d$ and $\mu>0$, a Nesterov smoothing approximation of the function
$f\colon\R^d\rightarrow\R$ defined by
\begin{equation*}
f(x):=\|x-a\|,\quad x\in\R^d,
\end{equation*}
admits the smooth DC representation
\begin{equation*}
f_\mu(x):=\dfrac{1}{2\mu}\|x-a\|^2-\dfrac{\mu}{2}\bigg[d\bigg(\dfrac{x-a}{\mu};\B\bigg)\bigg]^2.
\end{equation*}
Furthermore, we have the relationships
\begin{equation*}
\nabla f_\mu(x)=P\bigg(\dfrac{x-a}{\mu};\Bbb B\bigg)\;\;\mbox{\rm and }\;
f_\mu(x)\le f(x)\le f_\mu(x)+\dfrac{\mu}{2},
\end{equation*}
where $\Bbb B\subset\R^d$ is the closed unit ball, and where $P$ stands for the Euclidean projection \eqref{proj}.
\end{Proposition}
Using Proposition~\ref{Prop3.1} allows us to approximate the objective function $\mathcal{F}$ in \eqref{main problem} by a smooth DC function $\mathcal{F}_\mu$ as $\mu>0$ defined as follows:
\begin{eqnarray*}
\begin{array}{ll}
\mathcal{F}_\mu(\mathbf U,\mathbf V)&:=\disp\dfrac{1}{2\mu}\sum_{i=1}^{k}\sum_{j=1}^{n}u_{i,j}^2\|a_{j}-v_{i}\|^2-\disp\dfrac{\mu}{2}\sum_{i=1}^{k}\sum_{j=1}^{n}u_{i,j}^2\bigg[d\bigg(\dfrac{a_{j}-v_{i}}{\mu};\B\bigg)\bigg]^2
\\&=\mathcal{G}_\mu(\mathbf U,\mathbf V)-\mathcal{H}_\mu(\mathbf U,\mathbf V),
\end{array}
\end{eqnarray*}
where $\mathcal{G}_\mu,\mathcal{H}_\mu\colon\R^{k\times n}\times\R^{k\times d}\to\R$ are given by
\begin{eqnarray*}
\begin{array}{ll}
&\mathcal{G}_\mu(\mathbf U,\mathbf V):=\disp\dfrac{1}{2\mu}\sum_{i=1}^{k}\sum_{j=1}^{n}u_{i,j}^2\|a_{j}-v_{i}\|^2,\\
&\mathcal{H}_\mu(\mathbf U,\mathbf V):=\disp\dfrac{\mu}{2}\sum_{i=1}^{k}\sum_{j=1}^{n}u_{i,j}^2\bigg[d\bigg(\dfrac{a_{j}-v_{i}}{\mu};\B\bigg)\bigg]^2 \nonumber.
\end{array}
\end{eqnarray*}
This leads us to the construction of the following family of smooth approximations of the main problem \eqref{main problem} defined by
\begin{eqnarray}\label{discrete}
\begin{array}{ll}
&\mbox{\rm minimize }\;\mathcal{F}_\mu(\mathbf U,\mathbf V):=\mathcal{G}_\mu(\mathbf U,\mathbf V)-\mathcal{H}_\mu(\mathbf U,\mathbf V)\;\;\mbox{\rm as }\;\mu>0\\
&\mbox{\rm subject to }\;\mathbf U\in\mathcal{U}=\Delta^{n}\cap\{0,1\}^{k \times n}\;\;\mbox{\rm and }\;\mathbf V\in\R^{k\times d},
\end{array}
\end{eqnarray}
where $\Delta^n$ is the the $n$th Cartesian degree of the $(k-1)$-simplex $\Delta:=\{y\in[0,1]^k\;|\;\sum_{i=1}^{k}y_{i}=1\}$, which is a subset of $\R^{k}$.

Observe that for each $\mu>0$ problem \eqref{discrete} is of {\em discrete optimization}, while our intention is to convert it to a family of problems of {\em continuous optimization} for which we are going to develop and implement a DCA-based algorithm in Section~\ref{TheAlgorithm}.

The rest of this section is devoted to deriving two results, which justify such a reduction. The first theorem allows us to verify the {\em existence} of optimal solutions to the constrained optimization problems that appear in this procedure. It is required for having well-posedness of the algorithm construction.

\begin{Theorem}\label{existence} Let $(\Bar{\mathbf U},\Bar{\mathbf V})$ be an optimal solution to problem \eqref{discrete}. Then for any $\mu>0$ we have $\Bar{\mathbf V}\in\mathcal{B}$, where $\mathcal{B}:=\prod_{i=1}^{k}B_{i}$ is the Cartesian product of the $k$ Euclidean balls $B_{i}$ centered at $0\in\R^d$ with radius $r:=\sqrt{\sum_{j=1}^{n}\|a_j\|^2}$ that contain the optimal centers $\bar v_i$ for each index $i=1,\ldots,k$.
\end{Theorem}
{\bf Proof.} We can clearly rewrite the objective function in \eqref{discrete} in the form
\begin{equation}\label{F}
\mathcal{F}_\mu(\mathbf U,\mathbf V)=\dfrac{1}{2\mu}\sum_{i=1}^{k}\sum_{j=1}^{n}u_{i,j}\|a_{j}-v_{i}\|^2-\dfrac{\mu}{2}\sum_{i=1}^{k}\sum_{j=1}^{n}u_{i,j}\bigg[d\bigg(\dfrac{a_{j}-v_{i}}{\mu};\B\bigg)\bigg]^2
\end{equation}
due to interchangeability between $u_{i,j}^2$ and $u_{i,j}$. Observe that $\mathcal{F}_\mu(\mathbf U,\mathbf V)$ is differentiable on $R^{k\times n}\times\R^{k\times d}$. Employing the classical Fermat rule in \eqref{discrete} with respect to $V$ gives us $\nabla_{\mathbf V}\mathcal{F}_\mu(\Bar{\mathbf U},\Bar{\mathbf V})=0$. To calculate this partial gradient, we need some clarification for the second term in \eqref{F}, which is differentiable as a whole while containing the nonsmooth distance function \eqref{dist}. The convexity of the distance function in the setting of \eqref{F} allows us to apply the subdifferential calculation of convex analysis (see, e.g., \cite[Theorem~2.39]{bmn}) and to combine it with an appropriate chain rule to handle the composition in \eqref{F}. Observe that the distance function square in \eqref{F} is the composition of the nondecreasing convex function $\ph(t):=t^2$ on $[0, \infty)$ and the distance function to the ball $\B$. Thus the chain rule from \cite[Corollary~2.62]{bmn} is applicable. Thus, we can show that $d^2(\cdot; \B)$ is differentiable with
\begin{equation}\label{sdsub}
\nabla d^2(x; \B)=2[x-P(x; \B)]\; \mbox{\rm for }x\in \R^d.
\end{equation}
Using \eqref{sdsub}, we consider the following two cases:

\textbf{Case~1:} $(a_{j}-\bar v_{i})/\mu\in\B$ for the fixed indices $i\in\{1,\ldots,k\}$ and $j\in\{1,\ldots,n\}$. Then
\begin{equation*}
\nabla d^2\bigg(\dfrac{a_{j}-\bar v_{i}}{\mu};\B\bigg)=\{0\},
\end{equation*}
which  gives us
\begin{equation*}
\disp\frac{\partial\mathcal{F}_\mu}{\partial v_i}(\Bar{\mathbf U},\Bar{\mathbf V})=\dfrac{1}{\mu}\sum_{j=1}^{n}\bar u_{i,j}(\bar v_{i}-a_{j}),\quad i=1,\ldots,k,
\end{equation*}
for the corresponding partial derivatives of $\mathcal{F}_\mu$.

{\bf Case~2:} $(a_{j}-\bar v_{i})/\mu\notin\B$ for the fixed indices $i\in\{1,\ldots,k\}$ and $j\in\{1,\ldots,n\}$. In this case we have
\begin{eqnarray*}
\disp\frac{\partial\mathcal{F}_\mu}{\partial v_i}(\Bar{\mathbf U},\Bar{\mathbf V})&=&\dfrac{1}{2\mu}\sum_{j=1}^{n}\bar u_{i,j}2(\bar v_{i}-a_{j})+\sum_{j=1}^{n}\bar u_{i,j}\bigg[\dfrac{a_{j}-\bar v_{i}}{\mu}-P\bigg(\dfrac{a_{j}-\bar v_{i}}{\mu};\B\bigg)\bigg]\\
&=&\dfrac{1}{\mu}\sum_{j=1}^{n}\bar u_{i,j}(\bar v_{i}-a_{j})+\disp\sum_{j=1}^{n} \bar u_{i,j}\bigg[\dfrac{a_{j}-\bar v_{i}}{\mu}-\bigg(\dfrac{a_{j}-\bar v_{i}}{\|a_{j}-\bar v_{i}\|}\bigg)\bigg]\\
&=&\disp\dfrac{1}{\|a_{j}-\bar v_{i}\|}\sum_{j=1}^{n}\bar u_{i,j}(\bar v_{i}-a_{j}).
\end{eqnarray*}
Thus in both cases above it follows from the stationary condition $\nabla_{\mathbf V}\mathcal{F}_\mu(\Bar{\mathbf U},\Bar{\mathbf V})=0$ that
\begin{eqnarray*}
\bar v_{i}=\dfrac{\sum_{j=1}^{n}\bar u_{i,j}a_{j}}{\sum_{j=1}^{n}\bar u_{i,j}}\;\;\mbox{\rm for all }\;i=1,\ldots,k,
\end{eqnarray*}
since we have $\sum\limits_{j=1}^{n}\bar u_{i,j}>0$ due to the nonemptiness of the clusters. Then the classical Cauchy-Schwarz inequality leads us to the estimates
\begin{equation*}
\|\bar v_{i}\|^2 \le\dfrac{\bigg(\sum_{j=1}^{n}\bar u_{i,j}a_{j}\bigg)^2}{\bigg(\sum_{j=1}^{n}\bar u_{i,j}\bigg)^2}\le\sum_{j=1}^{n}\|a_{j}\|^2:=r^2,
\end{equation*}
which therefore verify all the conclusions of this theorem. $\h$\vspace*{0.05in}

Our next step is to enclose each discrete optimization problem \eqref{discrete} into the corresponding one of {\em continuous optimization}. For the reader's convenience if no confusion arises, we keep the same notation ${\mathbf U}$ for {\em all} the $k\times n$-matrices without the discrete restrictions on their entries. Define now the function $\mathcal{P}\colon\R^{k\times n}\to\R$ by
\begin{equation*}
\mathcal{P}(\mathbf U):=\sum_{i=1}^{k}\sum_{j=1}^{n}u_{i,j}(1-u_{i,j})\;\;\mbox{\rm for all }\;{\mathbf U}\in\R^{k\times n}
\end{equation*}
and observe that this function is concave on $\R^{k\times n}$ with $\mathcal{P}(\mathbf U)\ge 0$ whenever $\mathbf U\in\Delta^n$. Furthermore, we have the representations
\begin{equation}\label{constr}
\mathcal{U}=\big\{\mathbf U\in\Delta^n\;\big|\;\mathcal{P}(\mathbf U)=0\big\}=\big\{\mathbf U\in\Delta^n\;\big|\;\mathcal{P}(\mathbf U)\le 0\big\}
\end{equation}
for the set of feasible $k\times n$-matrices $\mathcal{U}$ in the original problem \eqref{main problem}. Employing further the standard {\em penalty function} method allows us to eliminate the most involved constraint on ${\mathbf U}$ in \eqref{constr} given by the function ${\mathcal P}$. Taking the penalty parameter $\alpha>0$ sufficiently large and using the smoothing parameter $\mu>0$ sufficiently small, consider the following family of continuous optimization problems:
\begin{eqnarray}\label{discrete1}
\begin{array}{ll}
&\mbox{\rm minimize }\;\mathcal{F}_\mu(\mathbf U,\mathbf V)+\alpha\mathcal{P}(\mathbf U)=\mathcal{G}_\mu(\mathbf U,\mathbf V)-\mathcal{H}_\mu(\mathbf U,\mathbf V)+\alpha\mathcal{P}(\mathbf U)\\
&\mbox{\rm subject to }\;\mathbf U\in\Delta^n\;\;\mbox{\rm and }\;\mathbf V\in\mathcal{B}.
\end{array}
\end{eqnarray}
Observe that Theorem~\ref{existence} ensures the existence of feasible solutions to problem \eqref{discrete1} and hence optimal solutions to this problem by the Weierstrass theorem due to the continuity of the objective functions therein and the compactness of the constraints sets $\Delta^n$ and ${\mathcal B}$.

Let us introduce yet another parameter $\rho>0$ ensuring a DC representation of the objective function in \eqref{discrete1} as follows:
\begin{eqnarray*}
\begin{array}{ll}
\mathcal{F}_\mu(\mathbf U,\mathbf V)+\alpha\mathcal{P}(\mathbf U)&=\dfrac{\rho}{2}\|(\mathbf U,\mathbf V)\|^2-\bigg(\dfrac{\rho}{2}\|(\mathbf U,\mathbf V)\|^2-\mathcal{F}_\mu(\mathbf U,\mathbf V)-\alpha\mathcal{P}(\mathbf U)\bigg)\\\nonumber
&=\dfrac{\rho}{2}\|(\mathbf U,\mathbf V)\|^{2}-\bigg(\dfrac{\rho}{2}\|(\mathbf U,\mathbf V)\|^{2}-\mathcal{G}_\mu(\mathbf U,\mathbf V)+\mathcal{H}_\mu(\mathbf U,\mathbf V)-\alpha\mathcal{P}(\mathbf U)\bigg)\\
&=:\mathcal{G}(\mathbf U,\mathbf V)-\mathcal{H}(\mathbf U,\mathbf V),
\end{array}
\end{eqnarray*}
where the function $\mathcal{G}(\mathbf U,\mathbf V):=\dfrac{\rho}{2}\|(\mathbf U,\mathbf V)\|^2$ is obviously convex, and
\begin{equation*}
\mathcal{H}(\mathbf U,\mathbf V):=\dfrac{\rho}{2}\|(\mathbf U,\mathbf V)\|^2-\mathcal{G}_\mu(\mathbf U,\mathbf V)+\mathcal{H}_\mu(\mathbf U,\mathbf V)-\alpha\mathcal{P}(\mathbf U).
\end{equation*}
Since $\mathcal{H}_\mu(\mathbf U,\mathbf V)-\alpha\mathcal{P}(\mathbf U)$ is also convex as $\alpha>0$, we are going to show that for any given number $\mu>0$ it is possible to determine the values of the parameter $\rho>0$ such that the function $\dfrac{\rho}{2}\|(\mathbf U,\mathbf V)\|^2-\mathcal{G}_\mu(\mathbf U,\mathbf V)$ is convex under an appropriate choice of $\rho$. This would yield the convexity of $\mathcal{H}(\mathbf U,\mathbf V)$ and therefore would justify a desired representation of the objective function in \eqref{discrete1}. The following result gives us a precise meaning of this statement, which therefore verify the required reduction of \eqref{discrete1} to {\em DC continuous optimization}.

\begin{Theorem}\label{convexity} The function
\begin{eqnarray}\label{g1}
\mathcal{G}_1(\mathbf U,\mathbf V):=\dfrac{\rho}{2}\|(\mathbf U,\mathbf V)\|^2-\mathcal{G}_\mu(\mathbf U,\mathbf V)
\end{eqnarray}
is convex on $\Delta^n\times\mathcal{B}$ provided that
\begin{equation}\label{rho}
\rho\ge\dfrac{n}{2\mu}\bigg[\bigg(1+\dfrac{1}{n}{\xi}^2\bigg)+\sqrt{\bigg(1+\dfrac{1}{n}{\xi}^2\bigg)^2+\dfrac{12}{n}{\xi}^2}\bigg],
\end{equation}
where $\xi:=r+\underset{1\le j\le n}{\mbox{max}}\|a_{j}\|$ and $r:=\sqrt{\sum_{j=1}^{n}\|a_j\|^2}$.
\end{Theorem}
{\bf Proof.} Consider the function $\mathcal{G}_1(\mathbf U,\mathbf V)$ defined in \eqref{g1} for all $(\mathbf U,\mathbf V)\in\Delta^n\times\mathcal{B}$ and deduce by elementary transformations directly from its construction that
\begin{eqnarray*}
\mathcal{G}_1(\mathbf U,\mathbf V)&=&\dfrac{\rho}{2}\|(\mathbf U,\mathbf V)\|^2-\mathcal{G}_\mu(\mathbf U,\mathbf V)\nonumber\\
&=&\dfrac{\rho}{2}\|(\mathbf U,\mathbf V)\|^2 - \dfrac{1}{2\mu}\sum_{i=1}^{k}\sum_{j=1}^{n}u_{i,j}^2\|a_{j}-v_{i}\|^2\nonumber\\
&=&\dfrac{\rho}{2}\|\mathbf U\|^2+\dfrac{\rho}{2}\|\mathbf V\|^2-\dfrac{1}{2\mu}\sum_{i=1}^{k}\sum_{j=1}^{n}u_{i,j}^2\|a_{j}-v_{i}\|^2\nonumber\\
&=&\dfrac{\rho}{2}\sum_{i=1}^{k}\sum_{j=1}^{n}u_{i,j}^2+\dfrac{\rho}{2n}\sum_{i=1}^{k}\sum_{j=1}^{n}\|v_{i}\|^2-\dfrac{1}{2\mu}\sum_{i=1}^{k}\sum_{j=1}^{n}u_{i,j}^2\|a_{j}-v_{i}\|^2\nonumber\\
&=&\sum_{i=1}^{k}\sum_{j=1}^{n}\dfrac{\rho}{2}u_{i,j}^2+\dfrac{\rho}{2n}\|a_{j}-v_{i}\|^2+\dfrac{\rho}{n}\langle a_j,v_i\rangle-\dfrac{\rho}{2n}\|a_{j}\|^2-\dfrac{1}{2\mu}u_{i,j}^{2}\|a_{j}-v_{i}\|^2.
\end{eqnarray*}
Next we define the functions $\gamma_{i,j}\colon\R\times\R^d\rightarrow\R$ for all $i=1,\ldots,k$ and $j=1,\ldots,n$ by
\begin{equation}\label{gamma}
\gamma_{i,j}(u_{i,j},v_{i}):=\dfrac{\rho}{2}u_{i,j}^2+\dfrac{\rho}{2n}\|a_{j}-v_{i}\|^2-\dfrac{1}{2\mu}u_{i,j}^2\|a_{j}-v_{i}\|^2
\end{equation}
and show that each of these functions is convex on the set $\{u_{i,j}\in[0,1],\,\|v_{i}\|\le r\}$, where $r>0$ is taken from Theorem~\ref{existence}.

To proceed, consider the Hessian matrix of each function in \eqref{gamma} given by
\begin{equation*}
J_{\gamma_{i,j}}(u_{i,j},v_{i}):=
\begin{bmatrix}
\rho-\dfrac{1}{\mu}\|a_{j}-v_{i}\|^2&&-\dfrac{2}{\mu}u_{i,j}(v_{i}-a_{j})\\
-\dfrac{2}{\mu}u_{i,j}(v_{i}-a_{j})&&\dfrac{\rho}{n}-\dfrac{1}{\mu}u^2_{i,j}
\end{bmatrix}
\end{equation*}
and calculate its determinant $\det(J_{\gamma_{i,j}}(u_{i,j},v_{i}))$ by
\begin{eqnarray*}
\label{det}
\det(J_{\gamma_{i,j}}(u_{i,j},v_{i})) &:=& \bigg(\rho-\dfrac{1}{\mu}\|a_{j}-v_{i}\|^2\bigg)\bigg(\dfrac{\rho}{n}-\dfrac{1}{\mu}u^2_{i,j}\bigg)-\dfrac{4}{\mu^{2}}u^2_{i,j}(v_{i}-a_{j})^T(v_{i}-a_{j})\nonumber\\
&=&\dfrac{\rho^{2}}{n}-\rho\bigg(\dfrac{u^2_{i,j}}{\mu}+\dfrac{1}{n\mu }\|v_{i}-a_{j}\|^{2}\bigg)-\dfrac{3u^2_{i,j}}{\mu^2}\|v_{i}-a_{j}\|^2.
\end{eqnarray*}
It follows from the well-known second-order characterization of the convexity that the function $\gamma_{i,j}(u_{i,j},v_{i})$ is convex on $\{u_{i,j}\in[0,1],\,\|v_{i}\|\le r\}$ if $\det(J_{\gamma_{i,j}}(u_{i,j},v_{i}))\ge 0$. Using \cite[Theorem~1]{Tao} gives us the estimate
\begin{eqnarray*}
\det\big(J_{\gamma_{i,j}}(u_{i,j},v_{i})\big)&\ge&\dfrac{\rho^2}{n}-\rho\bigg(\dfrac{1}{\mu}+\dfrac{1}{n\mu}\|v_{i}-a_{j}\|^2\bigg)-\dfrac{3}{\mu^2}\|v_{i}-a_{j}\|^2.
\end{eqnarray*}
Then we get from the construction of ${\mathcal B}$ in Theorem~\ref{existence} that $0<\|v_{i}-a_{j}\|\le\|v_{i}\|+\|a_{j}\|\le r+\underset{1\le j\le n}{\mbox{max}}\|a_{j}\|=:\xi$, and therefore
\begin{equation}
\label{quadratic}
\det\big(J_{\gamma_{i,j}}(u_{i,j},v_{i})\big)\ge\dfrac{\rho^2}{n}-\dfrac{\rho}{\mu}\bigg(1+\dfrac{1}{n}\xi^2\bigg)-\dfrac{3}{\mu^2}\xi^2.
\end{equation}
It allows us to deduce from the aforementioned condition for the convexity of $\gamma_{i,j}(u_{i,j},v_{i})$ that we do have this convexity if $\rho$ satisfies the estimate \eqref{rho}. $\h$

\section{Design and Implementation of the Solution Algorithm}~\label{TheAlgorithm}
Based on the developments presented in the previous sections and using the established smooth DC structure of problem \eqref{discrete1} with the subsequent $\rho-$parameterization of the objective function therein as ${\mathcal G}({\mathbf U},{\mathbf V})-{\mathcal H}({\mathbf U},{\mathbf V})$, we are now ready to propose and implement a new algorithm for solving this problem involving both DCA-2 and Nesterov's smoothing.

To proceed, let us present the problem under consideration in the equivalent {\em unconstrained} format by using the infinite penalty via the indicator function:
\begin{eqnarray}\label{unconstprob}
\begin{array}{ll}
\mbox{\rm minimize }&\dfrac{\rho}{2}\|(\mathbf U,\mathbf V)\|^2-\mathcal{H}(\mathbf U,\mathbf V)+\delta_{\Delta\times\mathcal{B}}(\mathbf U,\mathbf V)\\
&\mbox{\rm subject to }\;(\mathbf U,\mathbf V)\in\R^{k\times n}\times\R^{k\times d},
\end{array}
\end{eqnarray}
where $\mathcal{B}$, $\Delta$, and $\rho$ are taken from Section~\ref{nest-smooth}.

We first explicitly compute the gradient of the convex  function $\mathcal{H}(\mathbf U,\mathbf V)$ in \eqref{unconstprob}. Denoting
\begin{eqnarray*}
\begin{array}{ll}
[\mathcal{Y},\mathcal{Z}]&:=\nabla \mathcal{H}(\mathbf U,\mathbf V)=\nabla\bigg(\dfrac{\rho}{2}\|(\mathbf U,\mathbf V)\|^2-\disp\dfrac{1}{2\mu}\disp\sum_{i=1}^{k}\sum_{j=1}^{n}u_{i,j}^2\|a_{j}-v_{i}\|^2\\
&+\disp\dfrac{\mu}{2}\sum_{i=1}^{k}\sum_{j=1}^{n}u_{i,j}^2\bigg[d\bigg(\dfrac{a_{j}-v_{i}}{\mu};\B\bigg)\bigg]^2-\alpha\disp\sum_{i=1}^{k}\sum_{j=1}^{n}u_{i,j}(1-u_{i,j})\bigg),
\end{array}
\end{eqnarray*}
we have $\mathcal{Y}=\nabla\mathcal{H}_{\mathbf U}(\mathbf U,\mathbf V)$ and $\mathcal{Z}=\nabla\mathcal{H}_{\mathbf V}(\mathbf U,\mathbf V)$. Thus for each $i=1,\ldots,k$ and $j=1,\ldots,n$ the $(j,i)-$entry of the matrix $\mathcal{Y}$ and the $i$th row of the matrix $\mathcal{Z}$ are
\begin{eqnarray*}
\begin{array}{ll}
\mathcal{Y}_{j,i}&:=\rho u_{i,j}-\dfrac{u_{i,j}}{\mu}\|a_{j}-v_{i}\|^2+\mu u_{i,j}\bigg[d\bigg(\dfrac{a_{j}-v_{i}}{\mu};\B\bigg)\bigg]^2+2\alpha u_{i,j}-\alpha,\\
\mathcal{Z}_{i}&:=\disp\rho v_{i}-\dfrac{1}{\mu}\sum_{j=1}^{n}u_{i,j}^{2}(v_{i}-a_{j})-\sum_{j=1}^{n}u_{i,j}^2\bigg[\dfrac{a_{j}-v_{i}}{\mu}-P\bigg(\dfrac{a_{j}-v_{i}}{\mu};\B\bigg)\bigg],
\end{array}
\end{eqnarray*}
respectively. Let us now describe the proposed algorithm for solving the DC program \eqref{unconstprob} and hence the original problem \eqref{main problem} of multifacility location. The symbols $\mathcal{Y}^{l-1}_{[j,:]}$ and $\mathcal{Z}^{l-1}_{i}$ in this description represents the $j$th row of the matrix $\mathcal{Y}$ and the $i$th row of the matrix $\mathcal{Z}$ at the $l$th iteration, respectively. Accordingly we use the symbols $\mathbf U^{l}_{[:,j]}$ and $\mathbf V^{l}_{i}$. Recall also that the Frobenius norm of the matrices in this algorithm is defined in \eqref{frob}.\vspace*{0.05in}

{\bf Algorithm~3: Solving Multifacility Location Problems}.
\begin{center}
\begin{tabular}{| l |}
\hline
{\small INPUT}: $\mathbf X$ (the dataset), $\mathbf V^{0}$ (initial centers), ClusterNum (number of clusters), $\mu>0$,\\$\beta$ (scaling parameter) $>0$, $N\in\N$\\
{\small INITIALIZATION}: $\mathbf U^{0}$, $\epsilon >0$, $\mu_{f}$ (minimum threshold for $\mu$) $>0$, $\alpha>0$, $\rho>0$, \\
tol (tolerance parameter) $=1$\\
{\bf while} {tol $>\epsilon$ \mbox{and} $\mu >\mu_f$}\\
{\bf for} $l=1,2,\ldots,N$\\
For $1\le i\le k$ and $1\le j\le n$ compute\\\qquad \qquad$\mathcal{Y}^{l-1}_{j,i}:=\rho u^{l-1}_{i,j}-\dfrac{u^{l-1}_{i,j}}{\mu}\left\|a_{j}-v^{l-1}_{i}\right\|^{2}+\mu u^{l-1}_{i,j}\bigg[d\bigg(\dfrac{a_{j}-v^{l-1}_{i}}{\mu};\B\bigg)\bigg]^{2}+2\alpha u^{l-1}_{i,j}-\alpha,$\\
\\
\qquad\qquad$\mathcal{Z}^{l-1}_{i}:=\rho v^{l-1}_{i}-\dfrac{1}{\mu}\disp\sum_{j=1}^{n}(u^{l-1}_{i,j})^{2}(v^{l-1}_{i}-a_{j})-\sum_{j=1}^{n}(u^{l-1}_{i,j})^{2}\bigg[\dfrac{a_{j}-v^{l-1}_{i}}{\mu}-P\bigg(\dfrac{a_{j}-v^{l-1}_{i}}{\mu};\B\bigg)\bigg].$ \\
For $1\le i\le k$ and $1\le j\le n$ compute\\
\qquad\qquad$\mathbf U^{l}_{[:,j]}:=P\bigg(\dfrac{\mathcal{Y}^{l-1}_{[j,:]}}{\rho};\Delta\bigg)$,\\
\qquad\qquad$\mathbf V^{l}_{i}:=P\bigg(\dfrac{\mathcal{Z}^{l-1}_{i}}{\rho};B_{i}\bigg)=
\begin{cases}
\dfrac{\mathcal{Z}^{l-1}_{i}}{\rho} &\mbox{if }\;|\mathcal{Z}^{l-1}_{i}||\le\rho r,\\
\dfrac{r\mathcal{Z}^{l-1}_{i}}{||\mathcal{Z}^{l-1}_{i}||}&\mbox{if }\;||\mathcal{Z}^{l-1}_{i}||>\rho r.
\end{cases}
$\\
{\bf end for}\\
{\small UPDATE:}\\
\qquad\qquad $\text{tol}:=\left\|[\mathbf U^{l},\mathbf V^{l}]-[\mathbf U^{l-1},\mathbf V^{l-1}]\right\|_{F}$\\
\qquad\qquad $\mu:=\beta\mu$.\\
{\bf end while}\\
{\small OUTPUT}: $[\mathbf U^{N},\mathbf V^{N}]$.\\
\hline
\end{tabular}
\end{center}

Next we employ Algorithm~3 to solve several multifacility location problems of some practical meaning. By trial and error we verify that the values chosen for $\mu$ determine the performance of the algorithm for each data set. It can be seen that very small values of the smoothing parameter $\mu$ may prevent the algorithm from clustering, and thus we gradually decrease these values. This is done via multiplying $\mu$ by some number $0<\beta<1$ and stopping when $\mu<\mu_f$. Note also that in the implementation of our algorithm we use the standard approach of choosing $\mathbf U^{0}$ by computing the distance between the point in question and each group center $\mathbf V^{0}$ and then by classifying this point to be in the group whose center is the closest to it by assigning the value of $1$, while otherwise we assign the value of $0$.\vspace*{0.05in}

Let us now present several numerical examples, where we compute the optimal centers by using Algorithm~3 via MATLAB calculations. Fix in what follows the values of $\mu=0.5$, $\beta=0.85$, $\epsilon=10^{-6}$, $\mu_{f} =10^{-6}$, $\alpha=30$, and $\rho=30$ unless otherwise stated. The objective function is the {\em total distance} from the centers to the assigned data point. Note that this choice of the objective function seems to be natural from practical aspects in, e.g., airline and other transportation industries, where the goal is to reach the destination via the best possible route available. This reflects minimizing the transportation cost.

In the following examples we implement the standard $k$-means algorithm in MATLAB using the in-built function kmeans().
\begin{Example}\label{Eg4.1}{\rm
Let us consider a data set with $14$ entries in $\R^2$ given by
\[
\mathbf X:=
\begin{bmatrix}
0  & 2 & 7  & 2 & 3 & 6 & 5 & 8 & 8  & 9 & 1 &7 & 0 &0 \\
3 & 2 & 1  & 4 & 3 & 2 & 3  & 1 & 3 & 2 & 1 & 4 & 4&1
\end{bmatrix}^{T}
\]
with the initial data defined by
\begin{eqnarray*}
\begin{array}{ll}
&\mathbf V^0\ := \left[
\begin{array}{lr}
7.1429 & 2.2857 \\
1.1429 & 2.5714
\end{array}\right]\mbox{is obtained from the $k$-means algorithm; see Table~\ref{Table1}, }\\
&\mbox{ClusterNum}:=2.
\end{array}
\end{eqnarray*}
Employing Algorithm~3, we obtain the optimal centers as depicted in Table~\ref{Table1} and Figure~\ref{Figure 4.1}.

\begin{table}[H]
\caption{Comparison between Algorithm~3 and $k$-means}
\centering
\begin{tabular}{l l l}
\hline\hline
\textbf{Method} & \hspace{10pt} \textbf{Optimal Center ($\mathbf V^{N}$)} & \hspace{10pt} \textbf{Cost Function} \\
\hline
$k$-means & \hspace{20pt} $\left[
\begin{array}{lr}
7.1429 & 2.2857 \\
1.1429 & 2.5714
\end{array}\right]$ & \hspace{30pt}  22.1637\\
\\
Algorithm~3 & \hspace{20pt}
$\left[
\begin{array}{lr}
7.2220 & 2.1802 \\
1.1886 & 2.5069
\end{array}\right]$ & \hspace{30pt} 22.1352\\
\hline
\end{tabular}
\label{Table1}
\end{table}

\begin{figure}[H]
\centering
\includegraphics[scale=0.5]{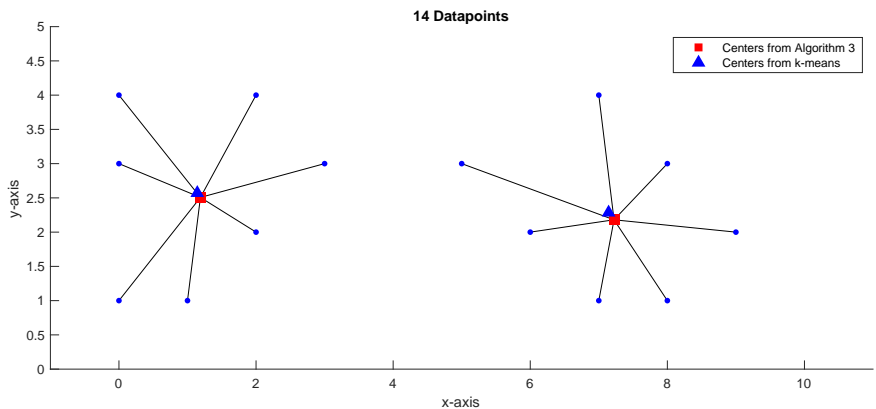}
\caption{\emph{MFLP with $14$ demand points and $2$ centers.}}
\label{Figure 4.1}
\end{figure}

Table~\ref{Table1} shows that the proposed Algorithm~3 is marginally better for the given data in comparison to the classical $k$-means approach in terms of the objective function.}
\end{Example}

\begin{Example}\label{balltests}{\rm In this example we test our algorithm on a dataset $\mathbf X$ containing 10 distinct points on each boundary of 4 balls of radius $r=0.3$ centered at (2,2), (4,2), (4,4), and (2,4). We generate the points as follows:
\begin{equation*}
\Big\{C_i+{r_i}\Big(\mbox{cos}\Big(\dfrac{j\pi}{5}\Big),\mbox{sin}\Big(\dfrac{j\pi}{5}\Big)\Big)\;\Big|\;\;i=1,\ldots,4;\;j=1,\ldots,10\Big\},
\end{equation*}
where $C_i$ and $r_i$ are the center and radius of each ball respectively; see \cite{NamGiles}. Typically the centroids are the centers of the balls. Choosing a random point from the boundary of each ball for the initial centers, this algorithm converges to the optimal solution
\[\mathbf V^{N}:=\begin{bmatrix}
2.0000&2.0000\\
4.0000&2.0000\\
2.0000&4.0000\\
4.0000&4.0000
\end{bmatrix}.\]
Its visualization is shown in Figure~\ref{Figure 6.2}.
\begin{figure}[H]
\centering
\includegraphics[scale=0.5]{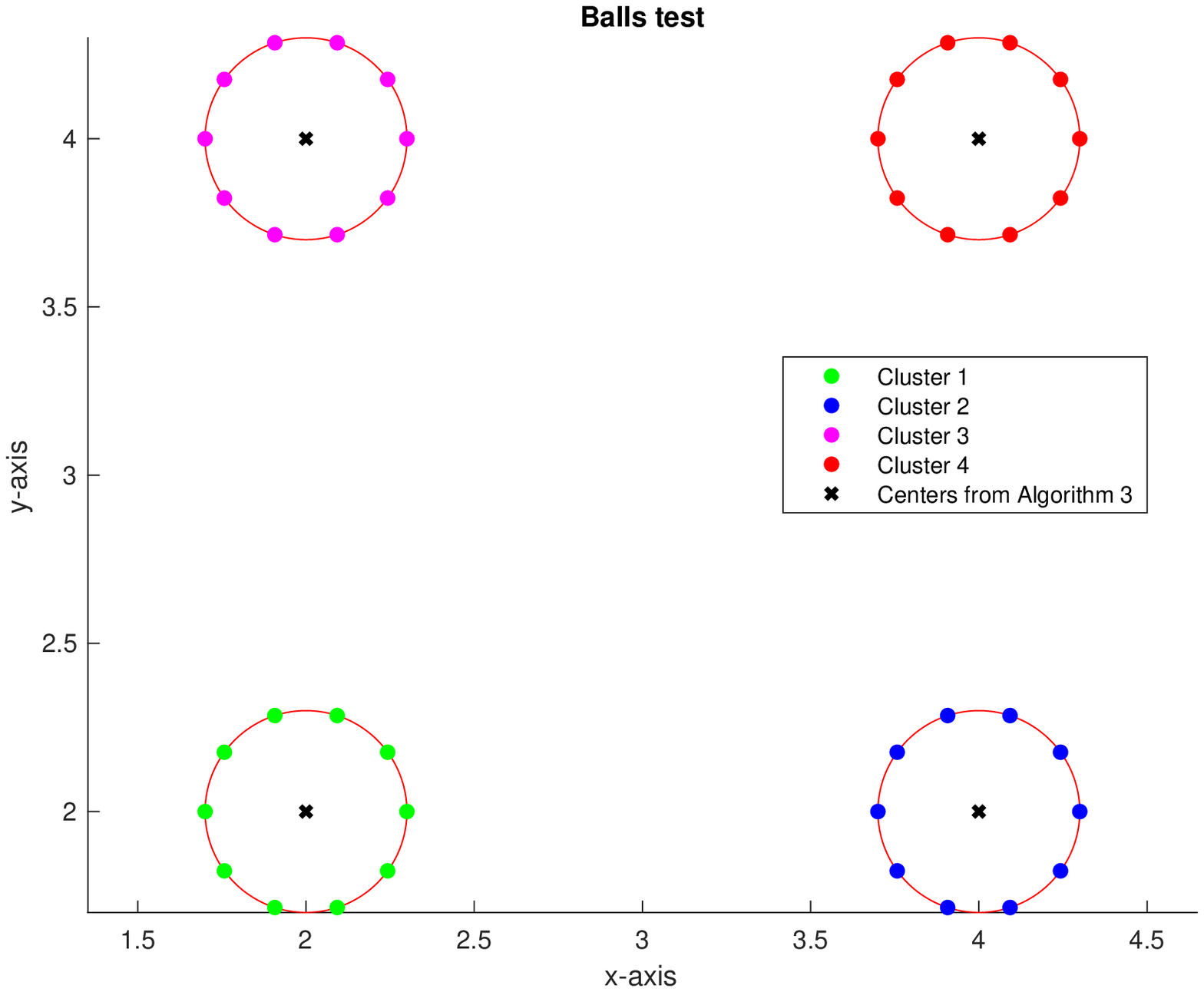}
\caption{\emph{MFLP with $40$ demand points and $4$ centers.}}
\label{Figure 6.2}
\end{figure}}
\end{Example}

Note that a drawback in employing the random approach to choose the initial cluster $\mathbf V^{0}$ in Example~5.2 is the need of having prior knowledge about the data. Typically it may not be plausible to extract such an information from large unpredictable real life datasets.

In the next example we choose the initial cluster by the process of {\em random selection} and see its effect on the optimal centers. Then the results obtained in this way by Algorithm~3 are compared with those computed by the $k$-means approach.

\begin{Example}\label{Eg4.2}{\rm
Let $\mathbf X$ be $200$ standard normally distributed random datapoints in $\R^2$, and let the initial data be given by
\begin{eqnarray*}
\begin{array}{ll}
&\mathbf V^0:=\mbox{randomly permuting and selecting $2$ rows of $\mathbf X$},\\
&\mbox{ClusterNum}:=2.
\end{array}
\end{eqnarray*}
We obtain the optimal centers as outlined in Table~\ref{TableAlg}.

\begin{table}[H]
\caption{Comparison between Algorithm~3 and $k$-means}
\centering
\begin{tabular}{l l l}
\hline\hline
\textbf{Method} & \hspace{10pt} \textbf{Optimal Center ($\mathbf V^N$)} & \hspace{10pt} \textbf{Cost Function} \\
\hline
$k$-means & \hspace{10pt} $\left[
\begin{array}{lr}
2.1016 & 1.2320 \\
-1.3060 & -1.0047
\end{array}\right]$ & \hspace{20pt} 403.3966\\
\\
Algorithm~3 & \hspace{10pt}
$\left[
\begin{array}{lr}
1.4902 & 0.7406\\
-1.3464 & -1.0716
\end{array}\right]$ & \hspace{20pt} 401.7506\\
\hline
\end{tabular}
\label{TableAlg}
\end{table}
Observe from Table~\ref{TableAlg} that the proposed Algorithm~3 is better for the given data in comparison to the standard $k$-means approach.
In addition, our approach gives a better approximation for the optimal center as shown in Figure~\ref{Figure 4.2}.

\begin{figure}[H]
\centering
\includegraphics[scale=0.5]{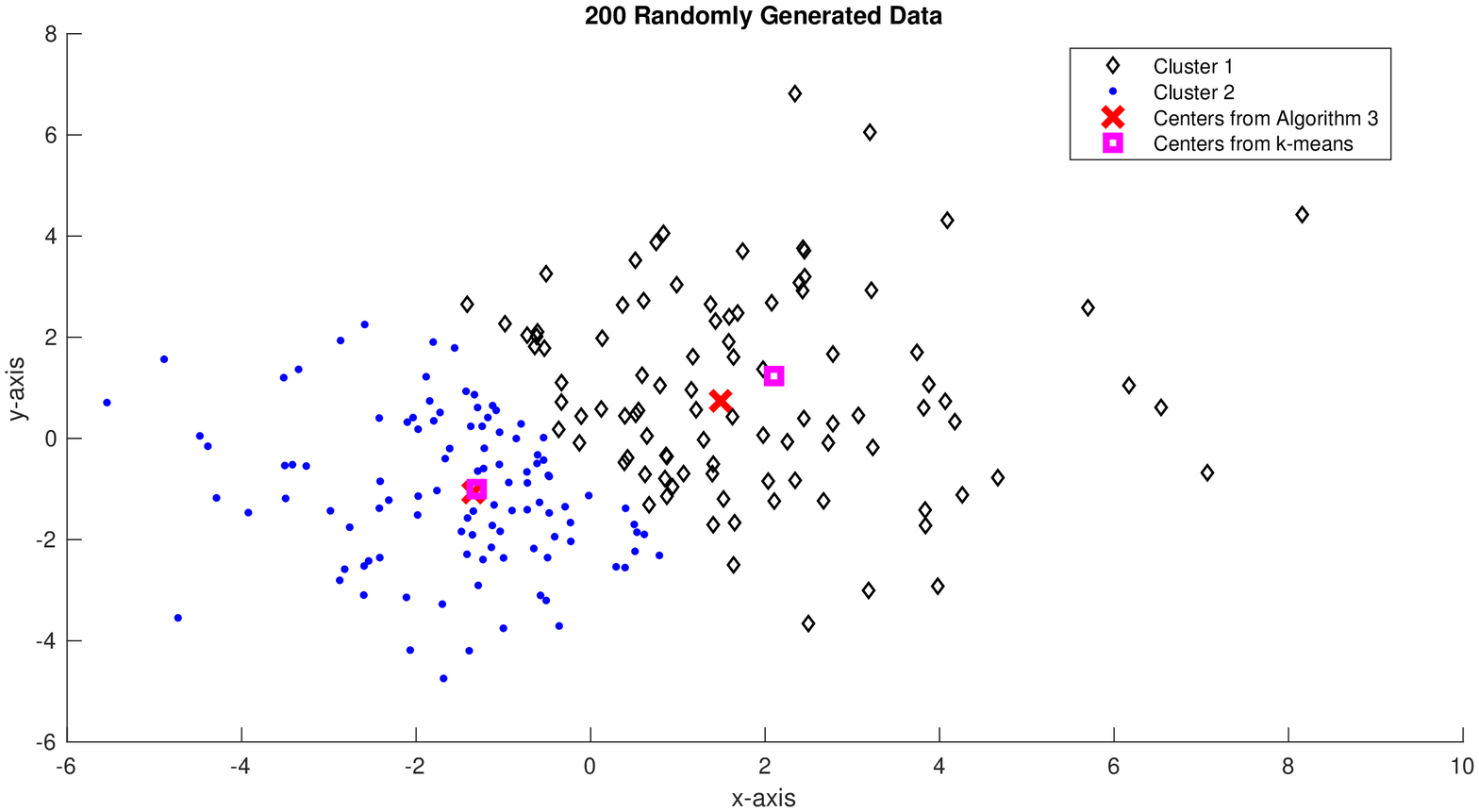}
\caption{\emph{MFLP with $200$ demand points and $2$ centers.}}
\label{Figure 4.2}
\end{figure}}
\end{Example}

Note that a real-life data may not be as efficiently clustered as in Example~\ref{Eg4.2}. Thus a suitable selection of the initial cluster $\mathbf V^0$ is vital for the convergence of the DCA~based algorithms. In the next Example~\ref{Eg4.3} we select $\mathbf V^0$ in Algorithm~3 by using the standard $k$-means method. The results achieved by our Algorithm~3 are again compared with those obtained by using the $k$-means approach.

\begin{Example}
\label{Eg4.3}{\rm
Consider the dataset $\mathbf X$ consisting of the latitudes and longitudes of $50$ most populous cities in the USA\footnote{Available at https://en.wikipedia.org/wiki/List of United States cities by population} with
\begin{eqnarray*}
	\begin{array}{ll}
		&\mathbf V^0\ := \left[
		\begin{array}{lr}
     	-80.9222 & 37.9882\\
    	-97.8273 & 35.3241 \\
	    -118.3121 & 36.9535
		\end{array}\right] \mbox{is obtained from the $k$-means algorithm; see Table~\ref{IntCluster}},\\
		&\mbox{ClusterNum}:=3.
	\end{array}
\end{eqnarray*}

By using Algorithm~3 we obtain the following optimal centers as given in Table~\ref{IntCluster}.

\begin{table}[H]
\caption{Comparison between Algorithm~3 (combined with $k$-means) and standard $k$-means}
\centering
\begin{tabular}{l l l }
\hline\hline
\textbf{Method} & \hspace{7pt}\textbf{Optimal Center ($\mathbf V^N$)} &  \textbf{Cost Function}\\
\hline
$k$-means
&\hspace{1pt}
$\left[
\begin{array}{lr}
-80.9222 & 37.9882\\
-97.8273 & 35.3241 \\
-118.3121 & 36.9535
\end{array}\right]$&\hspace{10pt} 288.8348\\
\\

Algorithm~3 (combined with $k$-means)&
\hspace{1pt}
$\left[
\begin{array}{lr}
-81.0970 & 38.3092\\
-97.4138 & 35.3383 \\
-119.3112 & 36.5410
\end{array}\right]$&\hspace{10pt} 286.6523\\
\hline
\end{tabular}
\label{IntCluster}
\end{table}

We see that Algorithm~3 (combined with $k$-means) in which the initial cluster $\mathbf V^0$ is selected by using $k$-means method performs better in comparison to the standard $k$-means approach (Table~\ref{IntCluster}). Moreover, it gives us optimal centers as depicted in Figure~\ref{Figure 4.3}.
\begin{figure}[H]
\centering
\includegraphics[scale=0.55]{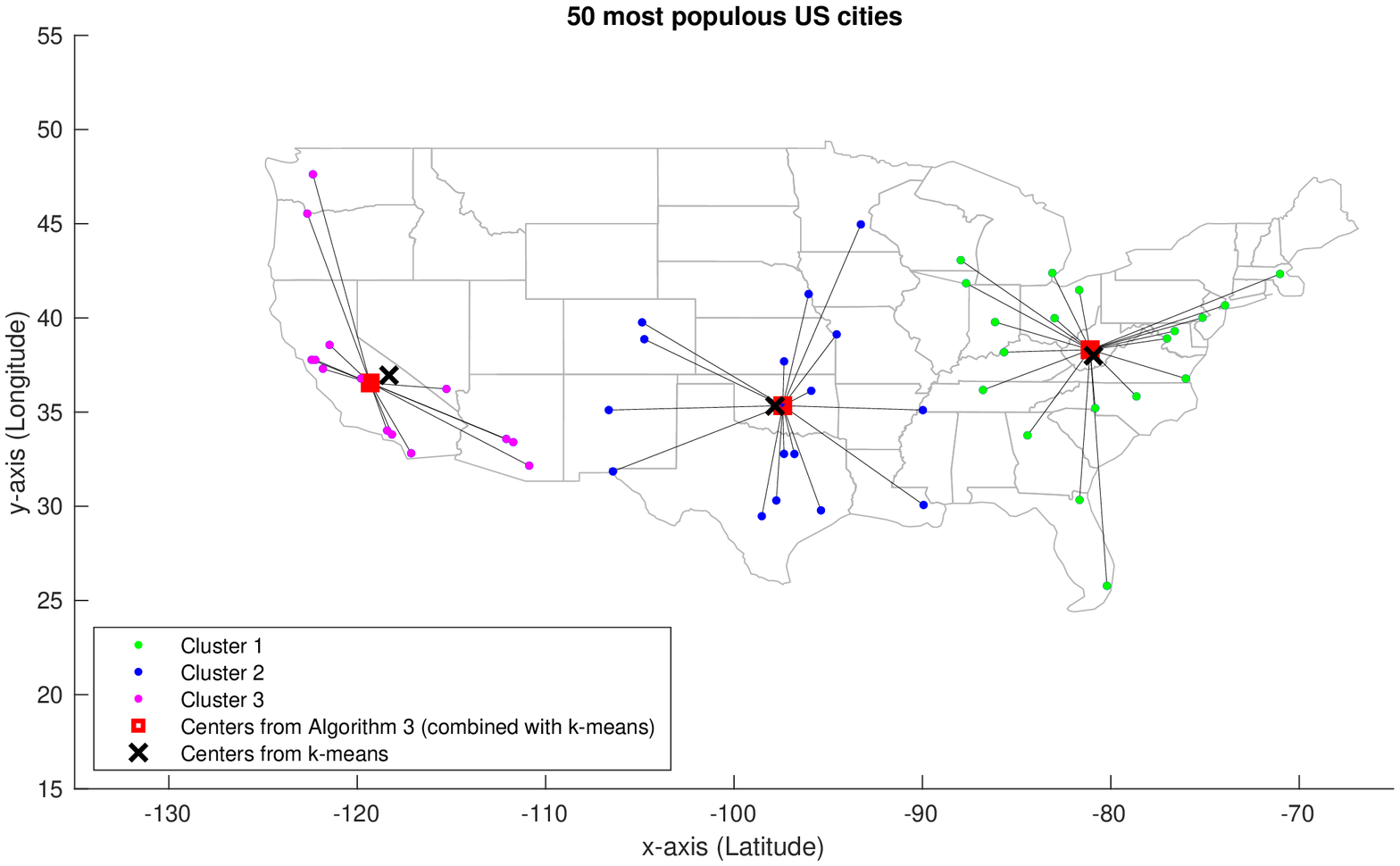}
\caption{\emph{MFLP with $50$ demand points and $3$ centers.}}
\label{Figure 4.3}
\end{figure}}
\end{Example}

In the next example we efficiently solve yet another multifacility location problem by using Algorithm~3.

\begin{Example} {\rm Consider the dataset $\mathbf X$ in $\R^2$ that consists of the latitudes and longitudes of $988$ US cities \cite{citydata} with \begin{eqnarray*}
			\begin{array}{ll}
				&\mathbf V^0\ := \left[
				\begin{array}{lr}
					-89.6747 & 41.1726\\
					-88.4834 & 30.2475 \\
					-118.4471 & 35.2843 \\
					-75.7890 & 40.0329 \\
					-114.1897 & 43.0798
				\end{array}\right] \mbox{is obtained from the $k$-means algorithm},\\
				&\mbox{ClusterNum}:=5.
			\end{array}
		\end{eqnarray*}

The optimal centers illustrated in Figure~\ref{Figure 4.5} are given by

\begin{equation*}
\mathbf V^N:=\left[
\begin{array}{lr}
-88.2248 & 40.7241\\
-88.0154 & 30.7041 \\
-120.0735 & 33.1941 \\
-74.5328 & 38.6996\\
-113.2341 & 42.3969
\end{array}\right].
\end{equation*}
The total transportation cost in this problem is 5089.5150.
\begin{figure}[H]
\centering
\includegraphics[scale=0.45]{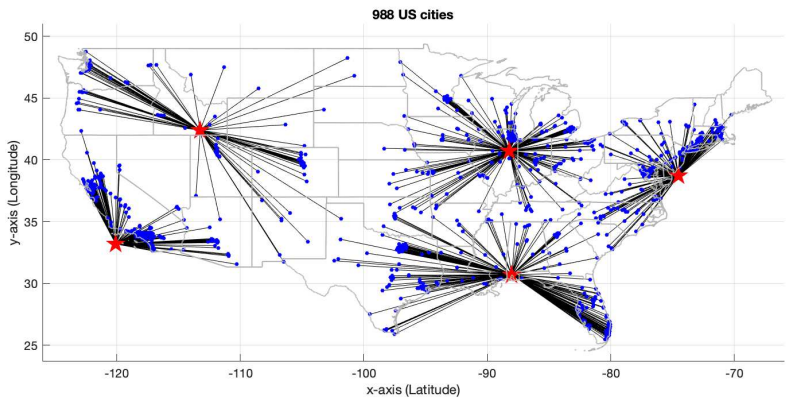}
\caption{\emph{MFLP with $988$ demand points and $5$ centers.}}
\label{Figure 4.5}
\end{figure}}
\end{Example}

In the last example presented in this section we efficiently solve a higher dimensional multifacility location problem by using Algorithm~3 and compare it's value of the cost function with the standard $k$-means algorithm.

\begin{Example} {\rm Let $\mathbf X$ in $\R^{13}$ be the \emph{wine dataset} from the UCI Machine Learning Repositiory \cite{wine} consisting of $178$ demand points. We apply Algorithm~3 with
		\begin{eqnarray*}
			\begin{array}{ll}
				&\mathbf V^0 \ \mbox{is obtained from the $k$-means algorithm},\\
				&\mbox{ClusterNum}:=3.
			\end{array}
		\end{eqnarray*}
	The total costs using Algorithm~3 and the $k$-means algorithm are obtained in Table~\ref{CostFunc} showing that former algorithm is better than the latter.
\begin{table}[H]
	\caption{Cost comparison between Algorithm~3 (combined with $k$-means) and $k$-means}
	\centering
	\begin{tabular}{ll}
		\hline\hline
		\textbf{Method} & 
		  \textbf{Cost Function}\\
		\hline
		$k$-means
	   &\hspace{10pt} 16556\\
		\\
		
		Algorithm~3 (combined with $k$-means)&\hspace{10pt} 16460\\
		\hline
	\end{tabular}
	\label{CostFunc}
\end{table}
}
\end{Example}

\section{Concluding Remarks}~\label{Conclusion}
In this paper we develop a new algorithm to solve a class of multifacility location problems. Its implementation exhibits a better approximation compared to the classical $k$-means approach. This is demonstrated by a series of examples dealing with two-dimensional problems with nonnegative weights. Thus the verification and implementation of the proposed algorithm for real-life multifacility location problems in higher dimensions with arbitrary weights is a central direction of our future work. In addition, refining the initial cluster selection and the stopping criterion is an important area to explore.


\begin{thebibliography}{99}

\bibitem{abt} L.T.H. An, M.T. Belghiti, and P.D. Tao, {\em A new efficient algorithm based on DC programming and DCA for clustering}, J. Glob. Optim. 37 (2007), pp.\ 593--608.

\bibitem{AnNam} N.T. An and N.M. Nam, {\em Convergence analysis of a proximal point algorithm for minimizing differences of functions}, Optim. 66 (2017), pp.\ 129--147.

\bibitem{TAN} L.T.H. An, H.V. Ngai, and P.D. Tao, {\em Convergence analysis of difference-of-convex algorithm with subanalytic data}, J. Optim. Theory Appl. 179 (2018), pp.\ 103--126.

\bibitem {ha} L.T.H.  An,  L.H. Minh, and P.D. Tao, {\em New and efficient DCA based algorithms for minimum sum-of-squares clustering}, Pattern Recogn. 47 (2014), pp.\ 388--401.

\bibitem {Tao} L.T.H. An and P.D. Tao, {\em Minimum sum-of-squares clustering by DC programming and DCA}, Int. Conf. Intell. Comp. (2009), pp.\ 327--340.
		
\bibitem{b} J. Brimberg, {\em The Fermat-Weber location problem revisited}, Math. Program. 71 (1995), pp.\ 71--76.

\bibitem{d} Z. Drezner, {\em On the convergence of the generalized Weiszfeld algorithm}, Ann. Oper. Res. 167 (2009), pp.\  327--336.

\bibitem{HM2015} T. Jahn, Y.S. Kupitz, H. Martini, and C. Richter, {\em Minsum location extended to gauges and to convex sets}, J. Optim. Theory Appl. 166 (2015), pp.\ 711--746.

\bibitem{HUL} J.-B. Hiriart-Urruty and C. Lemar\'{e}chal, {\em Fundamentals of Convex Analysis}, Springer, Berlin, 2001.

\bibitem{k} H.W. Kuhn, {\em A note on Fermat-Torricelli problem}, Math. Program. 4 (1973), pp.\ 98--107.

\bibitem{Ku-Ma} Y.S. Kupitz and H. Martini, {\em Geometric aspects of the generalized Fermat-Torricelli problem}, Bolyai Soc. Math. Stud. 6 (1997), pp.\ 55--128.

\bibitem{Martini} H. Martini, K.J. Swanepoel, and G. Weiss, {\em The Fermat-Torricelli problem in normed planes and spaces}, J. Optim. Theory Appl. 115 (2002), pp.\ 283--314.

\bibitem{m18} B.S. Mordukhovich, {\em Variational Analysis and Applications}, Springer, Cham, Switzerland, 2018.

\bibitem{n2} B.S. Mordukhovich and N.M. Nam, {\em Applications of variational analysis to a generalized Fermat-Torricelli problem}, J. Optim. Theory Appl. 148 (2011), pp.\ 431--454.

\bibitem{bmn} B.S. Mordukhovich and N.M. Nam, {\em An Easy Path to Convex Analysis and Applications}, Morgan \& Claypool Publishers, San Rafael, 2014.

\bibitem{nars} N.M. Nam, N.T. An, R.B. Rector, and J. Sun, {\em Nonsmooth algorithms and Nesterov's smoothing technique for generalized Fermat-Torricelli problems},
SIAM J. Optim. 24 (2014), pp.\ 1815--1839.

\bibitem {Nam2018} N.M. Nam, W. Geremew, S. Reynolds, and T. Tran, {\em Nesterov's smoothing technique and minimizing differences of convex functions for hierarchical clustering}, Optim. Lett. 12
(2018), pp.\ 455--473.

\bibitem{nh} N.M. Nam and N. Hoang, {\em A generalized Sylvester problem and a generalized Fermat-Torricelli problem}, J. Convex Anal. 20 (2013), pp.\ 669--687.
			
\bibitem{NamGiles} N.M. Nam, R.B. Rector, and D. Giles, {\em Minimizing differences of convex functions with applications to facility location and clustering}, J. Optim. Theory Appl. 173 (2017),
pp.\ 255--278.

\bibitem{n83} Yu. Nesterov, {\em A method for unconstrained convex minimization problem with the rate of convergence $O(1/k^2)$}, Soviet Math. Dokl. 269 (1983), pp.\ 543–-547.

\bibitem{n} Yu. Nesterov, {\em Smooth minimization of nonsmooth functions}, Math. Program. 103 (2005), pp.\ 127--152.

\bibitem{n18} Yu. Nesterov, {\em Lectures on Convex Optimization}, 2nd edition, Springer, Cham, Switzerland, 2018.

\bibitem{r} R.T. Rockafellar, {\em Convex Analysis}, Princeton University Press, Princeton, NJ, 1970.

\bibitem{TA1} P.D. Tao and L.T.H. An, {\em Convex analysis approach to D.C.\ programming: theory, algorithms and applications}, Acta Math. Vietnam. 22 (1997), pp.\ 289--355.

\bibitem{TA2} P.D. Tao and L.T.H. An, {\em A D.C.\ optimization algorithm for solving the trust-region subproblem}, {\em SIAM J. Optim.} 8 (1998), pp.\ 476--505.

\bibitem{citydata} United States Cities Database, {\em Simple Maps: Geographic Data Products}, 2017, \url{http://simplemaps.com/data/us-cities}.
	
\bibitem{w} E. Weiszfeld, {\em Sur le point pour lequel la somme des distances de $n$ points donn\'es est minimum}, T$\hat{\mbox{\rm o}}$hoku Math. J. 43 (1937), pp.\ 355--386.

\bibitem{wine}UCI Machine Learning Repository, \url{https://archive.ics.uci.edu/ml/datasets/wine}.
\end{thebibliography}
\end{document}